
\def\title#1{{\titlefont\noindent #1\bigskip}}

\def\author#1{{\largefont\noindent #1}\medskip}

\def\beginlinemode{\endmode
 \begingroup\obeylines\def\endmode{\par\endgroup}}
\let\endmode=\par

\newbox\theaddress
\def\address{\smallskip\beginlinemode\parindent 0in\getaddress}
{\obeylines
\gdef\getaddress #1 
 #2
 {#1\gdef\addressee{#2}%
   \global\setbox\theaddress=\vbox\bgroup\raggedright%
    \everypar{\hangindent2em}#2
   \def\endaddress{\egroup\endgroup \copy\theaddress \medskip}}}

\def\thanks#1{\footnote{}{\eightpoint #1}}

\long\def\Abstract#1{{\eightpoint\narrower\vskip\baselineskip\noindent
#1\smallskip}}

\def\skipfirstword#1 {}

\def\ir#1{\csname #1\endcsname}

\newdimen\currentht
\newbox\droppedletter
\newdimen\droppedletterwdth
\newdimen\drophtinpts
\newdimen\dropindent

\def\irrnSection#1#2{\edef\tttempcs{\ir{#2}}
\vskip6pt\penalty-3000
{\bf\noindent \expandafter\skipfirstword\tttempcs. #1}
\vskip6pt}

\def\irSubsection#1#2{\edef\tttempcs{\ir{#2}}
\vskip\baselineskip\penalty-3000
{\bf\noindent \expandafter\skipfirstword\tttempcs. #1}
\vskip6pt}

\def\irSubsubsection#1#2{\edef\tttempcs{\ir{#2}}
\vskip\baselineskip\penalty-3000
{\bf\noindent \expandafter\skipfirstword\tttempcs. #1}
\vskip6pt}

\def\References{\vbox to.25in{\vfil}\noindent{}{\bf References}
\vskip6pt\par}

\def\References{\vskip6pt\noindent{}{\bf References}
\vskip6pt\par}

\def\baselinebreak{\par \ifdim\lastskip<6pt
         \removelastskip\penalty-200\vskip6pt\fi}

\long\def\prclm#1#2#3{\baselinebreak
\noindent{\bf \csname #2\endcsname}:\enspace{\sl #3\par}\baselinebreak}

\def\Prf{\noindent{\bf Proof}: }

\def\rem#1#2{\baselinebreak\noindent{\bf \csname #2\endcsname}:\enspace }

\def\qed{{$\diamondsuit$}\vskip6pt}

\def\bibitem#1{\par\indent\llap{\rlap{\bf [#1]}\indent}\indent\hangindent
2\parindent\ignorespaces}

\long\def\eatit#1{}

\def\leftheadlinetext{}
\def\rightheadlinetext{}

\def\leftheadline{{\eightrm\folio\hfil \leftheadlinetext\hfil}}
\def\rightheadline{{\eightrm\hfil\rightheadlinetext\hfil\folio}}

\headline={\ifnum\pageno=1\hfil\else
\ifodd\pageno\rightheadline\else\leftheadline\fi\fi}

\def\tenpoint{\def\rm{\fam0\tenrm}
\textfont0=\tenrm \scriptfont0=\sevenrm \scriptscriptfont0=\fiverm
\textfont1=\teni \scriptfont1=\seveni \scriptscriptfont1=\fivei
\def\mit{\fam1} \def\oldstyle{\fam1\teni}
\textfont2=\tensy \scriptfont2=\sevensy \scriptscriptfont2=\fivesy
\def\cal{\fam2}
\textfont3=\tenex \scriptfont3=\tenex \scriptscriptfont3=\tenex
\def\it{\fam\itfam\tenit} 
\textfont\itfam=\tenit
\def\sl{\fam\slfam\tensl} 
\textfont\slfam=\tensl
\def\bf{\fam\bffam\tenbf} 
\textfont\bffam=\tenbf \scriptfont\bffam=\sevenbf
\scriptscriptfont\bffam=\fivebf
\def\tt{\fam\ttfam\tentt} 
\textfont\ttfam=\tentt
\normalbaselineskip=12pt
\setbox\strutbox=\hbox{\vrule height8.5pt depth3.5pt  width0pt}%
\normalbaselines\rm}

\def\eightpoint{\def\rm{\fam0\eightrm}%
\textfont0=\eightrm \scriptfont0=\sixrm \scriptscriptfont0=\fiverm
\textfont1=\eighti \scriptfont1=\sixi \scriptscriptfont1=\fivei
\def\mit{\fam1} \def\oldstyle{\fam1\eighti}%
\textfont2=\eightsy \scriptfont2=\sixsy \scriptscriptfont2=\fivesy
\def\cal{\fam2}%
\textfont3=\tenex \scriptfont3=\tenex \scriptscriptfont3=\tenex
\def\it{\fam\itfam\eightit} 
\textfont\itfam=\eightit
\def\sl{\fam\slfam\eightsl} 
\textfont\slfam=\eightsl
\def\bf{\fam\bffam\eightbf} 
\textfont\bffam=\eightbf \scriptfont\bffam=\sixbf
\scriptscriptfont\bffam=\fivebf
\def\tt{\fam\ttfam\eighttt} 
\textfont\ttfam=\eighttt
\normalbaselineskip=9pt%
\setbox\strutbox=\hbox{\vrule height7pt depth2pt  width0pt}%
\normalbaselines\rm}

\def\largefont{\def\rm{\fam0\largerm}
\textfont0=\largerm \scriptfont0=\largescriptrm \scriptscriptfont0=\tenrm
\textfont1=\largei \scriptfont1=\largescripti \scriptscriptfont1=\teni
\def\mit{\fam1} \def\oldstyle{\fam1\teni}
\textfont2=\largesy 
\def\cal{\fam2}
\def\it{\fam\itfam\largeit} 
\textfont\itfam=\largeit
\def\sl{\fam\slfam\largesl} 
\textfont\slfam=\largesl
\def\bf{\fam\bffam\largebf} 
\textfont\bffam=\largebf 
\def\tt{\fam\ttfam\largett} 
\textfont\ttfam=\largett
\normalbaselineskip=17.28pt
\setbox\strutbox=\hbox{\vrule height12.25pt depth5pt  width0pt}%
\normalbaselines\rm}

\def\titlefont{\def\rm{\fam0\titlerm}
\textfont0=\titlerm \scriptfont0=\largescriptrm \scriptscriptfont0=\tenrm
\textfont1=\titlei \scriptfont1=\largescripti \scriptscriptfont1=\teni
\def\mit{\fam1} \def\oldstyle{\fam1\teni}
\textfont2=\titlesy 
\def\cal{\fam2}
\def\it{\fam\itfam\titleit} 
\textfont\itfam=\titleit
\def\sl{\fam\slfam\titlesl} 
\textfont\slfam=\titlesl
\def\bf{\fam\bffam\titlebf} 
\textfont\bffam=\titlebf 
\def\tt{\fam\ttfam\titlett} 
\textfont\ttfam=\titlett
\normalbaselineskip=24.8832pt
\setbox\strutbox=\hbox{\vrule height12.25pt depth5pt  width0pt}%
\normalbaselines\rm}

\nopagenumbers

\font\eightrm=cmr8
\font\eighti=cmmi8
\font\eightsy=cmsy8
\font\eightbf=cmbx8
\font\eighttt=cmtt8
\font\eightit=cmti8
\font\eightsl=cmsl8
\font\sixrm=cmr6
\font\sixi=cmmi6
\font\sixsy=cmsy6
\font\sixbf=cmbx6

\font\largerm=cmr12 at 17.28pt
\font\largei=cmmi12 at 17.28pt
\font\largescriptrm=cmr12 at 14.4pt
\font\largescripti=cmmi12 at 14.4pt
\font\largesy=cmsy10 at 17.28pt
\font\largebf=cmbx12 at 17.28pt
\font\largett=cmtt12 at 17.28pt
\font\largeit=cmti12 at 17.28pt
\font\largesl=cmsl12 at 17.28pt

\font\titlerm=cmr12 at 24.8832pt
\font\titlei=cmmi12 at 24.8832pt
\font\titlesy=cmsy10 at 24.8832pt
\font\titlebf=cmbx12 at 24.8832pt
\font\titlett=cmtt12 at 24.8832pt
\font\titleit=cmti12 at 24.8832pt
\font\titlesl=cmsl12 at 24.8832pt

\tenpoint



\def\manyby{\hbox to.75in{\hrulefill}}
\hsize 6.5in 
\vsize 9.2in

\tolerance 3000
\hbadness 3000

\def\CMa{CM1}
\def\CMb{CM2}
\def\evain{E}
\def\vanc{H1}
\def\ravello{H2}
\def\fatpts{H3}
\def\sqrpts{HHF}
\def\Hi{Hi}
\def\nag{N1}
\def\nagB{N2}
\def\roe{R}
\def\Xu{X}

\def\C#1{\hbox{$\cal #1$}}

\def\pr#1{\hbox{{\bf P}${}^{#1}$}}
\def\leftheadlinetext{Brian Harbourne}
\def\rightheadlinetext{Nagata's Conjecture}

\title{On Nagata's Conjecture}

\author{Brian Harbourne}

\address
Department of Mathematics and Statistics
University of Nebraska-Lincoln
Lincoln, NE 68588-0323
email: bharbourne1@unl.edu
WEB: http://www.math.unl.edu/$\sim$bharbour/
\smallskip
September 14, 1999\endaddress
\vskip-\baselineskip

\thanks{\vskip -6pt
\noindent This work benefitted from a National Science Foundation grant.
In addition, I would like to thank
Professors Alan Peterson and Gary Meisters for
helpful discussions, and Professor Joaquim Ro\'e
for some helpful comments.
\smallskip
\noindent 1991 {\it Mathematics Subject Classification. } 
Primary 14N99, 
14J26. 
Secondary 14Q10. 
\smallskip
\noindent {\it Key words and phrases. }  Multiple points, plane curves,
fat points, Nagata's conjecture.\smallskip}

\vskip\baselineskip
\Abstract{Abstract: This paper gives an improved lower bound on
the degrees $d$ such that for general points
$p_1,\ldots,p_n\in\pr2$ and $m>0$ there is a plane curve of 
degree $d$ vanishing at each $p_i$ with multiplicity at least $m$.}
\vskip\baselineskip

\irrnSection{Introduction}{intro}
In this paper we work over an arbitrary algebraically closed field.
For positive integers $m$ and $n$, define $d(m,n)$ to be the least 
integer $d$ such that for general points $p_1,\ldots,p_n\in\pr2$
there is a curve of degree $d$ vanishing at each point $p_i$ with
multiplicity at least $m$. For $n\ge 10$, Nagata [\ir{nag}] conjectures that 
$d(m,n)>m\sqrt{n}$, and proves this when $n>9$ is a square.
(For $n\le 9$, applying methods of [\ir{nagB}] 
it can be shown that $d(m,n)=\lceil c_nm\rceil$,
where $c_n=1,1,1.5,2,2,12/5,21/8,48/17$ and $3$ 
for $n=1,\ldots,9$, resp. Recall for any real number $c$
that $\lfloor c\rfloor$ is the greatest integer less 
than or equal to $c$ and $\lceil c \rceil$ is the least 
integer greater than or equal to $c$;
in particular, $\lfloor c\rfloor\le c\le \lceil c \rceil$.)

Clearly, if $n'\le n$, then $d(m,n')\le d(m,n)$,
so we see from Nagata's result above that 
$d(m,n)\ge d(m,\lfloor \sqrt{n}\rfloor^2) 
> m\lfloor \sqrt{n}\rfloor$ for $n\ge 16$.
In fact, it is not hard to show directly for all $n\ge 1$
the slightly weaker
inequality $d(m,n)\ge m\lfloor \sqrt{n}\rfloor$;
similar reasoning shows $d(m,n)\ge mn/\lceil \sqrt{n}\rceil$ 
as well (see \ir{easybnd}). In certain ranges of $n$,
however, Ro\'e [\ir{roe}] has recently given a better bound: for
$n\ge3$ he shows that $d(m,n)\ge mr(n)$, where Ro\'e's constant 
$r(n)$ is defined as $r(n)=(n-1)\Pi_{i=2}^{n-1}(1-i/(n-1+i^2))$.
Ro\'e applies a procedure, called unloading, to an astute
sequence of specializations, to derive an algorithm 
for computing a value $R(m,n)$ depending on $m$ and $n$.
It turns out on general principles that 
$d(m,n)\ge R(m,n)$; the bound $d(m,n)\ge mr(n)$ is obtained
by showing that $R(m,n)\ge mr(n)$.

Although it seems hard actually to prove that $R(m,n)>m\sqrt{n}$
for $m<\sqrt{n}$, examples suggest that this is at least typically 
true. Indeed, a direct check shows for $2\le m\le 100$ 
that $R(m,m^2)$ is, plus or minus at most 1,
equal to $m^2+m/10$. (In a personal communication, Prof. Ro\'e 
has told me that in fact $R(m,m^2)\ge m^2+\lfloor m/10\rfloor$ 
for $m$ up to 200.) These and other data
indicate that Ro\'e's result 
$d(m,n)\ge R(m,n)$ is the best general bound
currently known when $m$ is not too large compared to $n$. 
For comparison, [\ir{evain}]
proves Nagata's conjecture for values of $m$ up to about $\sqrt{n/2}$
and [\ir{CMa}] and [\ir{CMb}] determine an exact value of $d(m,n)$ when 
$m\le 12$ and $n\ge 10$. This exact value agrees with conjectures 
(see [\ir{vanc}, \ir{Hi}, \ir{ravello}, \ir{CMa}, \ir{CMb}, 
\ir{sqrpts}]) which imply for $n\ge 10$
that $d(m,n)$ should be the least positive integer $d$ such that 
$d^2+3d+2-nm^2-mn>0$. When $n$ is an even square
and $m\ge(\sqrt{n}-2)/4$, this $d$ 
is precisely $m\sqrt{n}+(\sqrt{n}-2)/2$
(see [\ir{sqrpts}]), which unfortunately 
tends to be somewhat larger than $R(m,n)$. (In fact, 
the current paper resulted from this author's interest in whether 
Ro\'e's algorithm might in some cases be used to justify 
$d(m,n)=m\sqrt{n}+(\sqrt{n}-2)/2$ for $n$ an even square, in which case
the results of [\ir{sqrpts}] would give a resolution of the ideal
defining the fat point subscheme $mp_1+\cdots+mp_n$.)

In this paper, using a single specialization inspired
by Ro\'e's, we obtain (see \ir{mainthm})
$$d(m,n)\ge \lceil m\lambda_n\rceil,$$ 
where $\lambda_n$ denotes $n\lfloor\sqrt{n}\rfloor/
\lceil\sqrt{n}\lfloor\sqrt{n}\rfloor\rceil$.

In \ir{comps} we show that this is an improvement on
the bounds previously known. In particular, 
we verify that:
\item{$\bullet$} $\lambda_n\ge \lfloor\sqrt{n}\rfloor$ for all $n\ge 1$,
with equality if and only if $n$ or $n-1$ is a square;
\item{$\bullet$} $\lambda_n\ge n/\lceil\sqrt{n}\rceil$ for all $n\ge 1$,
with equality if and only if $n$, $n+1$ or $n+2$ is a square; 
\item{$\bullet$} $\lambda_n>r(n)$ for all $n\ge 3$; and
\item{$\bullet$} for each $n\ge 3$, that 
$\hbox{lim}_{m\to\infty}(m\lambda_n-R(m,n))=\infty$.

In \ir{conj} we show for certain values of $n$ with $m$ not too
large, that our bound implies Nagata's conjecture.

\irrnSection{Background}{surf}
We refer the reader to [\ir{fatpts}] for justification
and amplification of the material in this section.
Given essentially distinct points $p_1,\ldots,p_n\in\pr2$
(meaning for $i=0,\ldots,n-1$ that $p_{i+1}\in X_{i}$
where $X_0=\pr2$ and $\pi_{i+1}:X_{i+1}\to X_{i}$ 
is the blow up of $p_{i+1}$), we will denote $X_n$
simply by $X$, with the morphism $\pi:X\to X_0$
being the composition $\pi_n\circ\ldots\circ\pi_1$ of the blow ups.
The inverse image of $p_i$ with respect to $\pi_i$ 
is a divisor on $X_i$; the class of the total transform to $X$ of this 
divisor will be denoted $e_i$. The class of the total transform
to $X$ of a line in $\pr2=X_0$ will be denoted $e_0$.
The divisor class group on $X$ is then freely generated by
the classes $e_i, i=0,\ldots,n$, with the intersection form
being defined by $e_i\cdot e_j=0$ for $i\ne j$, $e_0^2=1$ and
$e_i^2=-1$ for $i>0$.

Define $d_X(m,n)$ to be the least $t$
such that $h^0(X,\C O_X(te_0-m(e_1+\cdots+e_n)))>0$. Then
$d(m,n)$ is the maximum value of $d_X(m,n)$ over all
essentially distinct sets of $n$ points of \pr2. (By semicontinuity,
it follows that $d(m,n)=d_X(m,n)$ for a general set
of distinct points $p_1,\ldots,p_n$.)
To give a bound $d(m,n)\ge d$, 
it clearly suffices to find a $d$ and an $X$ 
for which we can check $d_X(m,n)\ge d$
(i.e., for which
$h^0(X,\C O_X((d-1)e_0-m(e_1+\cdots+e_n)))=0$).
This follows, for example, if $X$ has a numerically
effective (also called {\it nef}) divisor $C$ such that 
$C\cdot ((d-1)e_0-m(e_1+\cdots+e_n))<0$.
The following lemma, which is well-known,
is, as we show, easy to prove this way.
(The slightly stronger result
$d(m,n)> m\lfloor \sqrt{n} \rfloor$
which follows from [\ir{nag}] as mentioned above, requires a related
but somewhat more involved argument.)

\prclm{Lemma}{easybnd}{Let $m$ and $n$ be positive integers.
Then we have:
\item{(a)} $d(m,n)\ge m\lfloor \sqrt{n} \rfloor$, and
\item{(b)} $d(m,n)\ge mn/\lceil \sqrt{n} \rceil$.}

\Prf To prove (a), choose distinct points $p_1,\ldots,p_{r^2}$
of a smooth irreducible plane $r$-ic $C'$, with
$r=\lfloor \sqrt{n} \rfloor$. Let $X$ be the surface
obtained by blowing up \pr2 at $p_1,\ldots,p_{r^2}$
and let $C$ be the class of the proper transform to $X$ of $C'$.
Then $C$ (being reduced and irreducible with $C^2\ge 0$)
is numerically effective; i.e., by definition
$C\cdot F\ge 0$ for every class $F$ on $X$ with 
$h^0(X,\C O_X(F))>0$ (we will refer to such a
class $F$ as an {\it effective class}). 
In particular, $d(m,n)e_0-m(e_1+\cdots+e_{r^2})$ is effective
since $d(m,n)\ge d(m,r^2)$ and since
$d(m,r^2)e_0-m(e_1+\cdots+e_{r^2})$ is effective, so we have
$d(m,n)r\ge mr^2$, and hence $d(m,n)\ge mr=m\lfloor \sqrt{n} \rfloor$.

To prove (b), choose distinct points $p_1,\ldots,p_n$
of a smooth irreducible plane $r$-ic $C'$, where this time
$r=\lceil \sqrt{n} \rceil$ and $X$ is the surface
obtained by blowing up \pr2 at $p_1,\ldots,p_n$
and $C$ is the class of the proper transform to $X$ of $C'$.
Then reasoning as above gives
$d(m,n)r\ge mn$, and hence the result.\qed

\irrnSection{The Main Result}{main}
In this section, we use a special arrangement of
essentially distinct points, similar to what
is used in [\ir{roe}], to which we will apply
an argument analogous to that used in the proof of \ir{easybnd}. 

\prclm{Proposition}{mainprop}{Let $d, n$ and $r$ be positive
integers such that $(r/d)^2\ge n$ and $r\le n$. Then
$d(m,n)\ge mnd/r$.}

\Prf Let $C_1$ be a smooth plane curve of degree $d$.
Choose any point $p$ such that $p_1\in C_1\subset X_0=\pr2$.
Let $X_1$ be the blow up of $X_0$ at $p_1$, and let $C_2$
be the proper transform of $C_1$. Then choose $p_2$ to be
the point of $C_2$ infinitely near to $p_1$.
Continue in this way, iteratively obtaining essentially
distinct points $p_i$, $i=1,\ldots,r$, where, for $1<i\le r$,
$p_i$ is the point of $C_i$ infinitely near to $p_{i-1}$
with respect to the blowing up $\pi_{i-1}:X_{i-1}\to X_{i-2}$
of $p_{i-1}$, with $C_i$ being the proper transform of $C_{i-1}$
with respect to $\pi_{i-1}$.
(Thus $C_1$ and $p_1$ determine $p_i$ for $1<i\le r$.)

If $n>r$, choose additional points $p_{r+1},\ldots,p_n$ so that
again each point $p_i$ is infinitely near to $p_{i-1}$ for
$i\ge r+1$ but so that $p_{r+1}$ is not on the proper transform
of $C_r$ and none of $p_i$, $i\ge r+1$ is on the proper transform 
to $X_{i-1}$ of the exceptional locus
of the blow up morphism $X_{i-2}\to X_{i-3}$ (i.e., $p_i$ is chosen
so that $e_{i-1}-e_i$ but not $e_{i-2}-e_{i-1}-e_i$ 
is effective). As usual, we denote $X_n$ by $X$;
$C$ will denote the 
class of the proper transform of $C_1$ to $X$.

Then $C$ is the class of a smooth, 
irreducible curve, as is each of 
$e_1-e_2, \ldots,e_{n-1}-e_n$ and $e_n$. 
By hypothesis,
$d^2n\le r^2$ and $r\le n$ and hence $d^2\le r$; 
using $d^2\le r$, it is not hard to verify
that $rde_0-d^2(e_1+\cdots+e_r)-(r^2-rd^2)e_{r+1}$
is the sum of $rC$ and
various nonnegative multiples of $e_i-e_{i+1}$ 
for  $1\le i\le r$ (here we have assumed
that $r<n$; we leave it to the reader to
consider the case that $r=n$).
But $d^2n\le r^2$ implies
$r^2-rd^2\ge (n-r)d^2$, hence
the class $D=rde_0-d^2(e_1+\cdots+e_n)$
is the sum of $rde_0-d^2(e_1+\cdots+e_r)-(r^2-rd^2)e_{r+1}$
and various nonnegative multiples of $e_n$
and of $e_i-e_{i+1}$ for $r+1\le i\le n-1$.
Thus $D$ is effective. 
But $D\cdot C=0$, $D\cdot (e_i-e_{i+1})=0$
for $i>0$ and $D\cdot e_n\ge0$, so $D$,
being a sum of effective 
classes which it meets nonnegatively, is nef.
Therefore, 
$d(m,n)rd-d^2mn=(d(m,n)e_0-
m(e_1+\cdots+e_n))\cdot D\ge 0$;
i.e., $d(m,n)\ge mnd/r$, as claimed. \qed

As a corollary we derive:

\prclm{Theorem}{mainthm}{Let $d$ and $m$ be positive
integers; then  $d(m,n)\ge 
\lceil m\lambda_n\rceil$.}

\Prf Apply \ir{mainprop} with
$d=\lfloor\sqrt{n}\rfloor$ and $r=
\lceil\sqrt{n}d\rceil$. We merely need to check
that $(r/d)^2\ge n$ and $r\le n$. 
Clearly, $(r/d)^2\ge (\sqrt{n}d/d)^2 = n$.
For the other inequality, we have 
$n=d^2+t$ where $0\le t \le 2d$, so 
$n\le d^2+t+(t/(2d))^2=(d+t/(2d))^2$,
hence $\sqrt{n}d\le(d+t/(2d))d=d^2+t/2\le n$;
therefore, $r=\lceil \sqrt{n}d \rceil \le n$ as required.\qed

\irrnSection{Comparisons}{comps}
We begin by comparing $\lambda_n$ with $\lfloor\sqrt{n}\rfloor$
and $n/\lceil\sqrt{n}\rceil$. We will use repeatedly
the easy fact that any integer $n\ge0$ can be (uniquely) 
written in the form $n=s^2+t$, where $s$ is a nonnegative 
integer and $0\le t\le 2s$ (indeed, $s$ is just
$\lfloor \sqrt{n}\rfloor$).

\prclm{Proposition}{Nagcomp}{Let $n=s^2+t$ where
$s>0$ and $0\le t\le 2s$ are integers; then:
\item{(a)} $\lambda_n\ge \lfloor\sqrt{n}\rfloor$,
with equality if and only if $t=0$ or $t=1$, and
\item{(b)} $\lambda_n\ge n/\lceil\sqrt{n}\rceil$,
with equality if and only if $t=0$, $t=2s-1$ or $t=2s$.}

\Prf First, note that $\lfloor\sqrt{n}\rfloor=s$ and that
$\sqrt{n}\le s+t/(2s)$. Thus $\lambda_n\ge 
ns/\lceil(s+t/(2s))s\rceil=ns/\lceil(s^2+t/2)\rceil$,
which is $(s^3+ts)/(s^2+t/2)$ if $t$ is even and
$(s^3+ts)/(s^2+(t+1)/2)$ if $t$ is odd. In particular,
we always have $\lambda_n\ge (s^3+ts)/(s^2+(t+1)/2)$.

To prove (a), we must show $\lambda_n\ge s$. If $t=0$,
it is easy to see that $\lambda_n=s$. If $t=1$,
then $s^2< s\sqrt{s^2+1}< s^2+1/2$
gives $\lambda_n=(s^3+s)/(s^2+1)=s$. For $t>1$,
it suffices to check that $(s^3+ts)/(s^2+(t+1)/2)>s$.
Similarly, (b) is clear when $t=0$. For
$t>0$, we have $n/\lceil\sqrt{n}\rceil=(s^2+t)/(s+1)$
so for $0<t<2s-1$ it suffices to check that $(s^2+t)s/(s^2+(t+1)/2)>
(s^2+t)/(s+1)$. We leave the cases $t=2s-1$ and $t=2s$
to the reader. \qed

We next want to compare $m\lambda_n$
with Ro\'e's bounds $mr(n)$ and $R(m,n)$. In order to 
deal with $R(m,n)$ it will be helpful to describe Ro\'e's
algorithm for computing it. 

We first develop some notation and terminology.
Let ${\bf w}=(m_1,\ldots,m_n)$ be a vector; 
then $p({\bf w})$ will denote the vector
obtained from ${\bf w}$ by putting the entries $m_j$ with $j>1$
into nonincreasing order. We will use ${\bf v}_i$ to denote the 
vector $(1,-1,\ldots,-1,0,\ldots,0)$,
where there are $i$ entries of $-1$.
Replacing every negative entry of a vector
by 0 we will call {\it rectification}. 
We will define $q_i({\bf w})$ to be ${\bf w}$, if,
with respect to the usual dot product,
${\bf w}\cdot {\bf v}_i\ge 0$; otherwise 
$q_i({\bf w})$ will be the rectification of 
${\bf w}+{\bf v}_i$.

Now let $n\ge 3$ be an integer; for each integer $i$ with $2\le i\le n-1$
we describe a routine $\C O_i$. Given 
a vector ${\bf w}=(m_1,\ldots,m_n)$
of nonnegative integers, let $g_i$ denote the composition $pq_ip$,
so $g_i({\bf w})=pq_ip({\bf w})$, and consider the sequence
$g_i({\bf w}), g_ig_i({\bf w}), \ldots$. It is easy to see that
eventually the sequence stabilizes at some vector 
which we will denote by $\C O_i({\bf w})$.

Ro\'e's algorithm, then, is to apply
consecutively the routines $\C O_2,\ldots, \C O_{n-1}$
to an initial input vector ${\bf w}=(m,\ldots,m)$; the value
$R(m,n)$ is the first entry of the vector
$\C O_{n-1}\cdots\C O_2({\bf w})$.

Since ${\bf w}=(m,\ldots,m)$ is of particular interest, 
in this case we will denote the first component
of $\C O_{i}\cdots \C O_{2}({\bf w})$ by $R_i(m,n)$ and set $R_1(m,n)=m$; 
thus $R_{n-1}(m,n)=R(m,n)$. The sum of the 2nd through 
$n$th components of $\C O_{i}\cdots \C O_{2}({\bf w})$ will be denoted
by $S_i(m,n)$ and we set $S_1(m,n)=(n-1)m$. Suppressing $m$ and $n$
when no confusion will result, we may write $R_i$ or $S_i$
instead of $R_i(m,n)$ or $S_i(m,n)$.

Another description of the algorithm 
will be helpful. Given integers $n>1$ and $1\le c\le n-1$, 
we will use $v(a,b,c,n)$ to denote the vector
$(a,b,\ldots,b,b-1,\ldots,b-1)$, where there are 
$c$ of the $b$ entries and $n$ entries altogether.
For example, if ${\bf w}=(m,\ldots,m)$, then 
${\bf w}=v(m,m,n-1,n)$.

Now let ${\bf v}=v(m_1,b,c,n)$ with $b\ge 1$ and assume that 
${\bf v}\cdot {\bf v}_{i-1}\ge 0$.
Then $g_i({\bf v})={\bf v}$ if ${\bf v}\cdot{\bf v}_i\ge0$
but if ${\bf v}\cdot{\bf v}_i<0$ then $g_i({\bf v})$
is $v(m_1+1,b',c',n)$, where $b'=b$ if $i<c$,
in which case $c'=c-i$, and $b'=b-1$ if $i\ge c$,
in which case $c'=n-i+c$. (Because of possible rectification,
the case that $b=1$ is a bit tricky, but since
${\bf v}\cdot{\bf v}_{i-1}\ge 0$, if $b=1$ and
${\bf v}\cdot{\bf v}_{i}<0$, then $i\le c$, so no
negative entries are ever involved.) 

In Ro\'e's algorithm, we continue to apply $g_i$ until
the dot product with ${\bf v}_i$ becomes nonnegative.
If starting with ${\bf v}$, $t$ is the least
number of such applications required for 
the dot product with ${\bf v}_i$ to become nonnegative,
then (denoting by $S$ the sum $(n-1)(b-1)+c$ 
of all components of ${\bf v}$ but the first)
the result of applying $g_i$ for $j\le t$ times is
$v(m_1+j, \lceil(S-ji)/(n-1)\rceil, \rho_j,n)$,
where $\rho_j=(S-ji)-(n-1)\lfloor(S-ji)/(n-1)\rfloor$ is the remainder 
when $S-ji$ is divided by $n-1$.
Looking at ${\bf v}_i\cdot v(m_1+j, \lceil(S-ji)/(n-1)\rceil, \rho_j,n)$ 
we see that $t$ is the least integer $j$ such that
$i\lfloor(S-ji)/(n-1)\rfloor+\hbox{min}(\rho_j,i)
\le m_1+j
\le i\lfloor(S-(j-1)i)/(n-1)\rfloor+\hbox{min}(\rho_{j-1},i)$.

In particular, applying the above remarks to
$v(R_{i-1},b,c,n)=\C O_{i-1}\cdots\C O_2({\bf w})$
and $v(R_i,b',c',n)=\C O_{i}\cdots\C O_2({\bf w})$, 
where ${\bf w}=(m,\ldots,m)$, we have the following formulas:
\itemitem{(F1)} $R_i=R_{i-1}+t$ where $t$ is the least $j$ such that 
$i\lfloor(S_{i-1}-ji)/(n-1)\rfloor+\hbox{min}(\rho_j,i)
\le R_{i-1}+j
\le i\lfloor(S_{i-1}-(j-1)i)/(n-1)\rfloor+\hbox{min}(\rho_{j-1},i)$;
\itemitem{(F2)} $S_{i-1}=(b-1)(n-1)+c$; and
\itemitem{(F3)} $S_i=S_{i-1}-i(R_i-R_{i-1})$ or equivalently 
$R_{i}+S_{i}/i=R_{i-1}+S_{i-1}/i$
(which are the same as $i(R_i-R_{i-1})=S_{i-1}-S_i$, 
which holds since $R$ increases
by 1 for each decrease in $S$ by $i$).

We note that the value $mr(n)$ can be obtained by a similar 
but ``averaged'' procedure, which requires working 
over the rationals. In place of ${\bf v}_i$ we have
$\overline{{\bf v}}_i=(1,-i/(n-1),\ldots,-i/(n-1))$, and in 
place of $q_i$ we have $\overline{q}_i$, where
$\overline{q}_i({\bf w})$ is ${\bf w}$ if
${\bf w}\cdot \overline{{\bf v}}_i\ge 0$; otherwise 
$\overline{q}_i({\bf w})$ is the rectification of 
${\bf w}+t\overline{{\bf v}}_i$, where $t$ is chosen so that
$({\bf w}+t\overline{{\bf v}}_i)\cdot\overline{{\bf v}}_i=0$.
We define $\overline{g}_i $ to be $p\overline{q}_ip$
(we use $p$ simply for analogy;
because of the averaging, nothing important would
be affected if we did not use it), and we take 
$\overline{\C O}_i({\bf w})$ to be the vector at which
the sequence $\overline{g}_i({\bf w}), 
\overline{g}_i\overline{g}_i({\bf w}),\dots$ stabilizes.
(Note that $\overline{\C O}_i({\bf w})=\overline{g}_i({\bf w})$
if neither ${\bf w}$ nor $\overline{q}_i({\bf w})$ has a negative 
entry.)

Now let $m$ be a positive integer, 
let $r_1(n)=1$, let $mr_i(n)$ be the first entry
of ${\bf w}_i=\overline{\C O}_i\cdots\overline{\C O}_2({\bf w}_1)$,
where ${\bf w}_1$ is the $n$-vector $(m,\ldots,m)$, let
$s_1(n)=n-1$ and let $ms_i(n)$ be the sum of all of 
the entries but the first of ${\bf w}_i$.
It is not hard by induction to check that
\itemitem{(f0)} ${\bf w}_i\cdot {\bf v}_i=0$ for $2\le i\le n-1$ and hence 
${\bf w}_i\cdot {\bf v}_{i+1}<0$ for $1\le i<n-1$, and that
\itemitem{(f1)} $r_i(n)=r_{i-1}(n)(i/(i-1))(1-i/(n-1+i^2))$, and
\itemitem{(f2)} $s_i(n)=((n-1)/i)r_i(n)$, from which it follows that
\itemitem{(f3)} $r_{i}(n)+s_{i}(n)/i=r_{i-1}(n)+s_{i-1}(n)/i$ and that
\itemitem{(f4)} $r_i(n)=(i^2/(n-1+i^2))(r_{i-1}+s_{i-1}/i)$.

\noindent By (f1), of course, we have 
$r_i(n)=\Pi_{j=2}^i(j(1-j/(n-1+j^2))/(j-1))$
and hence $r(n)=r_{n-1}(n)$.

\prclm{Proposition}{Roedev}{Let $n\ge 3$ be an integer; then:
\item{(a)} $r(n)\le \sqrt{n-1}-\pi/8+1/\sqrt{n-1}$, and
\item{(b)} $R(m,n)\le mr(n)+2(n-1)$.}

\Prf On behalf of easier reading, we will in this proof use $k$
to denote $n-1$. By direct check, (a) holds for $2\le k\le 3$.
So assume $k\ge 4$. In any case, $r(n)=k\Pi_{i=2}^{k}(1-i/(k+i^2))$, and 
[\ir{roe}] shows that $(r(n))^2=k\Pi_{i=1}^{k-1}(1-(i/(k+i^2))^2)$.
But $\hbox{log}(1-x)<-x$ holds for $0<x<1$, so we have
$\sum_{i=1}^{k-1}\hbox{log}(1-(i/(i^2+k))^2)
\le -\sum_{i=1}^{k-1}(i/(i^2+k))^2
= -\sum_{i\ge1}(i/(i^2+k))^2 +\sum_{i\ge k}(i/(i^2+k))^2$.
However, $\sum_{i\ge k}(i/(i^2+k))^2\le \int_{i\ge k-1}(x/(x^2+k))^2dx
\le \int_{i\ge k-1}x^{-2}dx = 1/(k-1)$, and from
[\ir{roe}] we see that $\sum_{i\ge1}(i/(i^2+k))^2=
\pi(-\pi+(\hbox{sinh}(2\sqrt{k}\pi))/(2\sqrt{k}))/(4\hbox{sinh}^2(\sqrt{k}\pi))
= -\pi^2/(4\hbox{sinh}^2(\sqrt{k}\pi)) +  
(\pi/(8\sqrt{k}))(\hbox{sinh}(2\sqrt{k}\pi))/\hbox{sinh}^2(\sqrt{k}\pi)$.
But $(\pi/(8\sqrt{k}))(\hbox{sinh}(2\sqrt{k}\pi))/\hbox{sinh}^2(\sqrt{k}\pi)
=(\pi/(4\sqrt{k}))(1+\hbox{exp}(-2\sqrt{k}\pi))/(1-\hbox{exp}(-2\sqrt{k}\pi))
\ge (\pi/(4\sqrt{k}))(1+\hbox{exp}(-2\sqrt{k}\pi))^2$, so
$-\sum_{i=1}^{k-1}(i/(i^2+k))^2\le 1/(k-1) +\pi^2/(4\hbox{sinh}^2(\sqrt{k}\pi)) - 
(\pi/(4\sqrt{k}))(1+\hbox{exp}(-2\sqrt{k}\pi))^2
\le 1/(k-1) +\pi^2/(4\hbox{sinh}^2(\sqrt{k}\pi)) - (\pi/(4\sqrt{k}))
\le 1.1/(k-1)-(\pi/(4\sqrt{k}))$,
where $\pi^2/(4\hbox{sinh}^2(\sqrt{k}\pi))<0.1/(k-1)$
follows from $(10\pi^2/4)(k-1)<\hbox{sinh}^2(\sqrt{k}\pi)$, 
which itself is easy to check (look at a graph first). 
Since $1.1/(k-1)-\pi/(4\sqrt{k})$ is negative for $k\ge 4$, 
the Taylor series for $\hbox{exp}(1.1/(k-1)-\pi/(4\sqrt{k}))$
is alternating so
$(r(n))^2\le k\hbox{exp}(1.1/(k-1)-\pi/(4\sqrt{k}))
\le k(1+1.1/(k-1)-\pi/(4\sqrt{k})+(1/2)(1.1/(k-1)-\pi/(4\sqrt{k}))^2)$,
but $1.1/(k-1)+(1/2)(1.1/(k-1)-\pi/(4\sqrt{k}))^2<
1.1/(k-1)+(1/2)(\pi/(4\sqrt{k}))^2<2/k$ so 
$(r(n))^2\le k(1+2/k-\pi/(4\sqrt{k}))$,
hence $r(n)\le \sqrt{k(1+2/k-\pi/(4\sqrt{k}))}
\le \sqrt{k}(1+1/k-\pi/(8\sqrt{k}))=\sqrt{n-1}+1/\sqrt{n-1}-\pi/8$.

Now consider (b). We begin by showing
$R_i\le (i^2R_{i-1}+iS_{i-1}+i^2+ik)/(i^2+k)$. 
Let $t=R_i-R_{i-1}$. By (F1), we have $R_{i-1}+t
\le i\lfloor(S_{i-1}-(t-1)i)/k\rfloor+\hbox{min}(\rho_{t-1},i)$
and hence $R_{i-1}+t
\le i\lfloor(S_{i-1}-(t-1)i)/k\rfloor+i$.
But $i\lfloor(S_{i-1}-(t-1)i)/k\rfloor+i\le 
i(S_{i-1}-ti)/k+(i^2+ki)/k$,
so solving for $t$ gives 
$t\le (iS_{i-1}-kR_{i-1}+i^2+ki)/(k+i^2)$
and therefore $R_i=R_{i-1}+t\le 
(iS_{i-1}+i^2R_{i-1}+i^2+ki)/(k+i^2)$, as claimed.

Now (given $m$ and $n$, and suppressing the $n$
notationally) it will be sufficient to prove
by induction for each $i$ 
that $R_i\le mr_i+2k$ and $R_i+S_i/i\le mr_i+s_i/i+2k$.
Note that $R_1=m\le m+2k=mr_1+2k$, and  
$R_1+S_1/1=nm\le nm+2k=mr_1+ms_1/1+2k$.
So assume that $R_{i-1}\le mr_{i-1}+2k$ and 
$R_{i-1}+S_{i-1}/(i-1)\le mr_{i-1}+ms_{i-1}/(i-1)+2k$
hold for some $i\ge2$. 

Since $R_{i-1}\le mr_{i-1}+2k$ and
$R_{i-1}+S_{i-1}/(i-1) \le mr_{i-1}+ms_{i-1}/(i-1)+2k$,
then $R_{i-1}+S_{i-1}/i \le mr_{i-1}+ms_{i-1}/i+2k$
must also hold, and using (F3) and (f3) we therefore have
$R_i+S_i/i=R_{i-1}+S_{i-1}/i
\le mr_{i-1}+ms_{i-1}/i+2k=mr_{i}+ms_{i}/i+2k$,
as required.

Since $R_{i-1}+S_{i-1}/i \le mr_{i-1}+ms_{i-1}/i+2k$ 
and since $R_i\le 
(iS_{i-1}+i^2R_{i-1}+i^2+ki)/(k+i^2)
=(i^2/(k+i^2))(R_{i-1}+S_{i-1}/i)+(i^2+ik)/(i^2+k)$,
the latter is at most
$(i^2/(k+i^2))(mr_{i-1}+ms_{i-1}/i)+(2ki^2+i^2+ik)/(i^2+k)$,
so by (f4),
this latter simplifies to $mr_i+(2ki^2+i^2+ik)/(i^2+k)$
which (taking $k$ to be $i$) 
is at most $mr_i+(2k^3+2k^2)/(k^2+k)=mr_i+2k$,
as we needed to show.
\qed

We now compare our bound with those of Ro\'e.

\prclm{Proposition}{Roecomp}{Let $n\ge 3$ be an integer; then:
\item{(a)} $\lambda_n>r(n)$, and
\item{(b)} $\hbox{lim}_{m\to\infty}m\lambda_n-R(m,n)=\infty$.
In particular, $\lceil m\lambda_n\rceil>R(m,n)$ for all
sufficiently large $m$.}

\Prf Let $s=\lfloor \sqrt{n}\rfloor$
and write $n=s^2+t$; thus $0\le t\le 2s$ and
$\lambda_n=(s^2+t)s/\lceil s\sqrt{n}\rceil\ge (s^2+t)s/\lceil s(s+t/(2s))\rceil
=(s^2+t)s/\lceil s^2+t/2\rceil\ge (s^2+t)s/(s^2+(t+1)/2)$.

For part (a), one first checks case by case that $r(n)<\lambda_n$ for 
$3\le n\le 48$, so we are reduced to the case that $n\ge 49$; i.e., $s\ge 7$. 
First assume $t=0$; then in fact $\lambda_n=\sqrt{n}$. But since
$n\ge 8$ we see $-\pi/8+1/\sqrt{n-1}<0$ so, by \ir{Roedev}(a),
$r(n)\le \sqrt{n-1}-\pi/8+1/\sqrt{n-1}<\sqrt{n}=\lambda_n$.

Hereafter we may assume that $t\ge 1$. Thus 
$r(n)\le \sqrt{n-1}-\pi/8+1/\sqrt{n-1}\le\sqrt{s^2+t-1}-\pi/8+1/s
\le s+t/(2s)-(1/(2s)+\pi/8-1/s)= s+t/(2s)+1/(2s)-\pi/8$,
and since $s\ge 4$ we see $1/(2s)-\pi/8<-1/4$.
Therefore $r(n)<\lambda_n$ follows if we show that 
$s+t/(2s)-1/4\le (s^2+t)s/(s^2+(t+1)/2)$, which simplifies to
$2s^2+t(t+1)\le s^3+(t+1)s/2$. But using $t\le 2s$ and $t\ge 1$,
respectively, we have $2s^2+t(t+1)\le 6s^2+2s$ and $s^3+s\le s^3+(t+1)s/2$,
so it is enough to show that $6s^2+2s\le s^3+s$, which is true for
$s\ge 7$.

Now (b) is clear: by \ir{Roedev}(b) we have $R(m,n)\le mr(n)+2n$,
and we have just checked that $r(n)<\lambda_n$. 
\qed

\irrnSection{Nagata's Conjecture}{conj}
Nagata's conjecture, that $d(m,n)>m\sqrt{n}$, has been verified 
for various small $m$: by [\ir{CMa},\ir{CMb}] for $m<13$ and
$n>9$, and for $m$ up to about $\sqrt{n/2}$ by [\ir{evain}],
and examples suggest that $R(m,n)>m\sqrt{n}$ for $m$
up to about $\sqrt{n}$. Our bound
$d(m,n)\ge \lceil m\lambda_n\rceil$ also implies 
$d(m,n)>m\sqrt{n}$
in certain situations, one such we give here.

\prclm{Theorem}{nc}{Let $n=s^2+s$, where $s\ge 3$ is an integer.
Then $d(m,n)>m\sqrt{n}$ holds if either 
\item{(a)} $s$ is even and $m<2s$, or
\item{(b)} $s$ is odd and $m<2s/3$.}

\Prf Since $d(m,n)\ge \lceil m\lambda_n\rceil$ and
$s+1/2>\sqrt{s^2+s}$, it suffices to
show that $\lceil m\lambda_n\rceil\ge m(s+1/2)$. This follows from
\ir{lnc}, for (b) using the fact that $2s/3<(2s^2+s+1)/(3s+1)$. \qed

\prclm{Lemma}{lnc}{Let $n=s^2+s$, where $s\ge 1$. 
\item{(a)} If $s$ is even, then 
$$\lceil m\lambda_n\rceil=\cases{ms+m/2, &if $0<m<4s+2$ and $m$ is even;\cr 
                                 ms+(m+1)/2, &if $0<m<2s$ and $m$ is odd.\cr}$$
\item{(b)}  If $s$ is odd, then 
$$\lceil m\lambda_n\rceil=\cases{ms+m/2, &if $0<m<(4s^2+2s+2)/(3s+1)$ and $m$ is even;\cr 
                                 ms+(m+1)/2, &if $0<m<(2s^2+s+1)/(3s+1)$ and $m$ is odd.\cr}$$}

\Prf For the moment we consider more generally
the case that $n=s^2+t$ for some $0\le t\le 2s$.
First say $t$ is even. Since $(s+t/(2s))\ge \sqrt{s^2+t}$, 
it follows that $s(s+t/(2s))=\lceil s\sqrt{s^2+t}\rceil$
since $(s^2+t/2)-s\sqrt{s^2+t}<1$.
Thus $m\lambda_n=ms(s^2+t)/(s^2+t/2)=ms(1+t/(2s^2+t))
=ms+mst/(2s^2+t)$. On the other hand, if $t$ is odd, then 
$s^2+(t+1)/2=\lceil s\sqrt{s^2+t}\rceil$ since
$0<(s^2+(t+1)/2)-s\sqrt{s^2+t}<1$.
This time we 
find $m\lambda_n=ms1+ms(t-1)/(2s^2+t+1))$.

(a) Using $t=s$ in our formula for $t$ even above, we 
have $\lceil m\lambda_n\rceil=ms+\lceil ms/(2s+1)\rceil$. If $m$ is even, then
$\lceil ms/(2s+1)\rceil=m/2$ if $m/2-ms/(2s+1)<1$, which holds for
$0<m<4s+2$, while if $m$ is odd, then
$\lceil ms/(2s+1)\rceil=(m+1)/2$ if $0<(m+1)/2-ms/(2s+1)<1$, which holds for
$0<m<2s+1$. 

The proof of (b) is similar.
\qed

\rem{Remark}{xu} In closing we mention
that, given a reduced and irreducible plane curve
of degree $d$ through $n$ general points of multiplicity $m$,
Xu [\ir{Xu}] shows that $d>m\sqrt{n}-1/(2\sqrt{n-1})$. \qed

\References

\bibitem{\CMa} Ciliberto, C. and Miranda, R. {\it Degenerations
of planar linear systems}, 
J. Reine Angew. Math. 501 (1998), 191-220.

\bibitem{\CMb} Ciliberto, C. and Miranda, R. {\it Linear systems
of plane curves with base points of equal multiplicity}, 
preprint (1998: \hbox{http://xxx.lanl.gov/abs/math/9804018}),
to appear, Trans AMS.

\bibitem{\evain} Evain, L. {\it Une minoration du degre des courbes plane
singularits imposes}, 
Preprint ENS Lyon 212 (1998), 1-17.

\bibitem{\vanc} Harbourne, B. {\it The geometry of rational surfaces and Hilbert
functions of points in the plane},
Can.\ Math.\ Soc.\ Conf.\ Proc.\ 6 
(1986), 95--111.

\bibitem{\ravello} Harbourne, B. {\it Points in Good Position in \pr 2}, in:
Zero-dimensional schemes, Proceedings of the
International Conference held in Ravello, Italy, June 8--13, 1992,
De Gruyter, 1994.

\bibitem{\fatpts} \manyby. {\it Free Resolutions of Fat Point 
Ideals on \pr2}, JPAA 125 (1998), 213--234.

\bibitem{\sqrpts} Harbourne, B., Holay, S. and Fitchett, S.
{\it Resolutions of Ideals of Uniform Fat Point Subschemes of \pr2},
preprint (1999: \hbox{http://xxx.lanl.gov/abs/math/9906130}).

\bibitem{\Hi} Hirschowitcz, A.
{\it Une conjecture pour la cohomologie 
des diviseurs sur les surfaces rationelles generiques},
Journ.\ Reine Angew.\ Math. 397
(1989), 208--213.

\bibitem{\nag} Nagata, M. {\it On the 14-th problem of Hilbert}, 
Amer.\ J.\ Math.\ 33 (1959), 766--772.

\bibitem{\nagB} Nagata, M. {\it On rational surfaces, II}, Mem.\ Coll.\ Sci.\ 
Univ.\ Kyoto, Ser.\ A Math.\ 33 (1960), 271--293.

\bibitem{\roe} Ro\'e, J. {\it On the existence of plane curves
with imposed multiple points}, 
preprint (1998: \hbox{http://xxx.lanl.gov/abs/math/9807066}).

\bibitem{\Xu} Xu, G.
{\it Curves in \pr2 and symplectic packings},
Math. Ann. 299 (1994), 609--613.

\bye